\documentclass[journal]{IEEEtran}
\usepackage{amsmath,epsfig}
\usepackage{inputenc}
\usepackage[T1]{fontenc}
\usepackage[french,english]{babel}
\usepackage{amssymb,amsthm,color,bbm}
\usepackage{dsfont}
\usepackage{graphicx}
\usepackage{subfig}
\usepackage{version}
\usepackage{ifthen}
\usepackage{url}
\usepackage{enumerate}
\usepackage{algorithm,algorithmic}

\newcommand{\paramstar}{m^{\star}}

\newcommand{\trans}{^{T}}

\newcommand{\bigR}{\mathbb{R}}

\newcommand{\argmax}{\operatornamewithlimits{argmax}}

\newcommand{\CE}[4][]
{
\ifthenelse{\equal{#1}{}}{\mathbb{E}^{#2}_{#3}\left[#4\right]}{\mathbb{E}^{#2}_{#3}\left[#4\middle | #1\right]}
}

\newcommand{\Rset}{\mathbb{R}}

\newcommand{\eqdef}{\ensuremath{\stackrel{\mathrm{def}}{=}}}

\renewcommand{\theenumi}{(\roman{enumi})}
\def\Xset{\mathbb{X}}
\def\Yset{\mathbb{Y}}
\def\Nset{\mathbb{N}}

\newcommand{\pscal}[2]{\left \langle #1, #2 \right\rangle}

\newcommand{\eqsp}{\;}

\newcommand{\rme}{\mathrm{e}}

\newcommand{\XinitIS}[2][]{\ifthenelse{\equal{#1}{}}{\ensuremath{\rho_{#2}}}{\ensuremath{\check{\rho}_{#2}}}}

\newcommand{\CExpparam}[4]{\mathbb{E}_{#1}^{#2}\left[#3\middle | #4\right]}

\newcommand{\filt}[2][]%
{
\ifthenelse{\equal{#1}{}}{\ensuremath{\phi_{#2}}}%
{\ifthenelse{\equal{#1}{hat}}{\ensuremath{\phi^{N}_{#2}}}
{\ifthenelse{\equal{#1}{tilde}}{\ensuremath{\tilde{\phi}^{N}_{#2}}}
{\ifthenelse{\equal{#1}{tar}}{\ensuremath{\phi^{N,\mathrm{t}}_{#2}}}
{\ifthenelse{\equal{#1}{aux}}{\ensuremath{\phi^{N,\mathrm{a}}_{#2}}}
}
}
}
}
}
\newcommand{\mcbf}[2][]{%
\ifthenelse{\equal{#1}{}}{\overline{\mathcal{F}}_{#2}}{\overline{\mathcal{F}}_{#2}^{(#1)}}%
}

\newcommand{\unfilt}[2][]%
{
\ifthenelse{\equal{#1}{}}{\ensuremath{\gamma_{#2}}}%
{\ifthenelse{\equal{#1}{hat}}{\ensuremath{\gamma^{N}_{#2}}}
{\ifthenelse{\equal{#1}{tilde}}{\ensuremath{\tilde{\gamma}^{N}_{#2}}}
{\ifthenelse{\equal{#1}{tar}}{\ensuremath{\gamma^{N,\mathrm{t}}_{#2}}}
{\ifthenelse{\equal{#1}{aux}}{\ensuremath{\gamma^{N,\mathrm{a}}_{#2}}}
}
}
}
}
}
\newcommand{\CPE}[3][]
{\ifthenelse{\equal{#1}{}}{\mathbb{E}\left[\left. #2 \, \right| #3 \right]}{\mathbb{E}_{#1}\left[\left. #2 \, \right| #3 \right]}}

\newcommand{\Xsigma}[1][]%
{%
\ifthenelse{\equal{#1}{}}{\ensuremath{\mathcal{B}(\Xset)}}{\ensuremath{\mathcal{B}(\Xset^{#1})}}
}
\newcommand{\sumwght}[2][]{%
\ifthenelse{\equal{#1}{}}{\ensuremath{\Omega_{#2}}}{\ensuremath{\Omega_{#2}^{(#1)}}}}

\newcommand{\adjfunc}[4][]
{\ifthenelse{\equal{#1}{}}{\ifthenelse{\equal{#4}{}}{\vartheta_{#2}}{\vartheta_{#2}(#4)}}
{\ifthenelse{\equal{#1}{smooth}}{\ifthenelse{\equal{#4}{}}{\tilde{\vartheta}_{#2}}{\tilde{\vartheta}_{#2}(#4)}}
{\ifthenelse{\equal{#1}{fully}}{\ifthenelse{\equal{#4}{}}{\vartheta^\star_{#2}}{\vartheta^\star_{#2}(#4)}}{\mathrm{erreur}}}}}
\newcommand{\chunk}[4][]%
{\ifthenelse{\equal{#1}{}}{\ensuremath{{#2}_{#3:#4}}}{\ensuremath{#2^#1}_{#3:#4}}
}
\newcommand{\kiss}[3][]
{\ifthenelse{\equal{#1}{}}{p_{#2}}
{\ifthenelse{\equal{#1}{fully}}{p^{\star}_{#2}}
{\ifthenelse{\equal{#1}{smooth}}{\tilde{r}_{#2}}{\mathrm{erreur}}}}}

\newcommand{\Kiss}[3][]
{\ifthenelse{\equal{#1}{}}{P_{#2}}
{\ifthenelse{\equal{#1}{fully}}{P^{\star}_{#2}}
{\ifthenelse{\equal{#1}{smooth}}{\tilde{R}_{#2}}{\mathrm{erreur}}}}}

\newcommand{\post}[3][]%
{
\ifthenelse{\equal{#1}{}}{\ensuremath{\phi_{#2|#3}}}%
{\ifthenelse{\equal{#1}{hat}}{\ensuremath{\phi^{\mathrm{FFBS},N}_{#2|#3}}}
{\ifthenelse{\equal{#1}{tilde}}{\ensuremath{\phi^{\mathrm{FFBSi},N}_{#2|#3}}}
{\ifthenelse{\equal{#1}{tar}}{\ensuremath{\phi^{N,\mathrm{t}}_{#2|#3}}}}
}
}
}

\newcounter{hypA}

\newtheorem{pr}{Proposition}

\begin{document}

\title{Simultaneous Localization and Mapping Problem in Wireless Sensor Networks}

\author{Thierry~Dumont\thanks{Thierry Dumont,  Laboratoire de Math\'ematiques, Universit\'e Paris-Sud}
\thanks{thierry.dumont@math.u-psud.fr}~and~Sylvain~Le~Corff\thanks{Sylvain Le Corff, LTCI, TELECOM ParisTech and CNRS}
\thanks{sylvain.lecorff@telecom-paristech.fr}
\thanks{This work is partially supported by ID Services, 22/24 rue Jean Rostand, 91400 Orsay, France.}}

\maketitle

\begin{abstract}
Mobile device localization in wireless sensor networks is a challenging task. It has already been addressed when the WiFI propagation maps of the access points are modeled deterministically. However, this procedure does not take into account the environmental dynamics and also assumes an offline human training calibration. In this paper, the maps are made of an average indoor propagation model  combined with a perturbation field which represents the influence of the environment. This perturbation field is embedded with a prior distribution. The device localization is dealt with using Sequential Monte Carlo methods and relies on the estimation of the propagation maps. This inference task is performed online, i.e. using the observations sequentially, with a recently proposed online Expectation Maximization based algorithm. The performance of the algorithm are illustrated through Monte Carlo experiments.
\end{abstract}

\begin{IEEEkeywords}
Simultaneous localization and Mapping, Indoor localization, Received Signal Strength Indicator, WiFi, Signal Propagation.
\end{IEEEkeywords}

\IEEEpeerreviewmaketitle

\section{Introduction}

\label{sec:intro}

Wireless sensor networks \cite{gaura:2010} generally consist of a data acquisition network and a data distribution network, monitored and controlled by a management center. These networks have many applications such as environmental monitoring  (\cite{LCAV-CONF-2008}) or target tracking (\cite{Lau:2009,Bahl:2000,Chen:2005}). In this paper, we consider a WiFi communication network made up of a mobile device (such as a hand held mobile computer or a smartphone), a server and WiFi access points (APs). We are interested in the estimation of the localization of the mobile device in the environment using the signal strength of the surrounding APs. The mobile device collects the power of the signals and sends the data to the server which uses them to build an estimator of the device's position. The key step to provide such an estimator is to understand the behaviour of the WiFi signal strength for different positions in the environment. However, predicting the propagation of WiFi signals in an indoor environment is challenging since they are subject to many perturbations (\textit{e.g.} shadowing, reflection...).

Two main techniques exist to approximate the WiFi signal propagation map of each AP: the first ones use deterministic models based on the localization and characteristics of the surrounding APs as well as the localization of the obstacles involved in the environment, see for instance \cite{Gorce:2007}. Other famous techniques are based on a previous hand made offline training phase in which a human operator performs a site survey by measuring the received signal strength indicator (RSSI) from different APs at some fixed sampled points, see \cite{Bahl:2000,Evennou:2006}. However, representing the indoor propagation map using a deterministic model is challenging since several obstacles cannot be taken into account. On the contrary, the site survey method allows to build an accurate estimation of the signal strength, but only for a finite number of points. Nevertheless, \cite{ferris:2006} provides a method to extend these measures to the entire map using Gaussian processes techniques.

In this paper, we propose an estimation method that does not require any calibration procedure. The propagation maps are estimated online (\textit{i.e.} without storing the observations) using the data sent by the mobile device. 
%Since only the measures gathered by the mobile device are considered in the approximation of the propagation maps, the time dependency of the positions accuracy can be minimized. 
Any modification in the way the WiFi signals propagate inside the environment (due to new obstacles for instance) affects the data sent by the mobile device. Then, while these changes deteriorate the accuracy of localization systems using fixed estimators of the propagation maps, our system learns these changes by taking them into account in the construction of our map estimators. Thus, as illustrated in Section \ref{sig:sec:truedata}, the accuracy of our localization method improves with time instead of degrading.

A semiparametric statistical model is used: the propagation maps are made of a parametric average indoor model in addition to a non parametric perturbation field. This model combines a prior knowledge on the signal propagation with random perturbations due to the obstacles. Based on the data collected by the mobile device, parameters and perturbation field estimators can be defined. We simultaneously provide an estimator for the device position. The procedure relies on an online Expectation-Maximization (EM) based algorithm for the estimation of the propagation maps and on particle filtering for the estimation of the device position.

The structure of this paper is the following. Section \ref{sec:framework} describes the model and defines the notations. Section \ref{sec:Online EM} presents the online EM algorithm and Section \ref{sec:algo} gives a general algorithm for online inference in our Simultaneous localization and Mapping (SLAM) problem. Section \ref{sec:simulations} illustrates this algorithm with numerical experiments.

\section{Model and assumption}
\label{sec:framework}

Let $\{X_{t}\}_{t\geq 1}$ be the cartesian coordinates of the mobile device in a two-dimensional compact space. This continuous environment is discretized into a finite grid map, denoted by $\mathcal{C}$, for purposes of numerical computation. It is assumed that $\{X_{t}\}_{t\geq 1}$ is a Markov chain taking values in $\mathcal{C}$ with initial distribution $\nu$ and Markov transition matrix given, for all $(x,x')\in \mathcal{C}^2$, by 
\begin{equation}
\label{eq:def:q}
q_{x,x^{\prime}}\propto\rme^{-\|x-x^{\prime}\|^{2}/a}\eqsp,
\end{equation}
where $a\in\Rset_{+}^{\star}$ depends on the average speed of the mobile and is assumed to be known  and $\|\cdot\|$ denotes the usual euclidean norm in $\Rset^{2}$. Let $B$ be the number of APs, $|\mathcal{C}|$ be the cardinality of $\mathcal{C}$ and $F^{\star}$ be the $B\times |\mathcal{C}|$ matrix where $F^{\star}_{j,x}$ is the $j$-th AP expected signal strength at position $x$. At each time step $t$, the mobile device measures and sends to the server the observation $Y_{t}$ taking values in $\Rset^{B}$. 
%For any $1\leq j\leq B$, $Y_{t,j}$ denotes the $j$-th component of $Y_{t}$ and represents the signal strength (in $dBm$) of the $j$-th AP.  
For all $t\geq 0$, the observation $Y_{t} $ is given by
\begin{equation}\label{eq:model}
Y_{t} \eqdef F^{\star}_{.,X_{t}} + \varepsilon_{t}\eqsp,
\end{equation}
where $F^{\star}_{.,X_{t}}\eqdef \{F^{\star}_{j,X_{t}}\}_{j=1}^{B}$ and where $\{\varepsilon_{t}\}_{t\geq 0}$ is a sequence of i.i.d Gaussian random vectors, independent from $\{X_t\}_{t\ge 1}$, with mean $0$ and covariance matrix $\Sigma\eqdef \sigma^{\star,2}I_{B}$ ($I_{B}$ is the identity matrix of size $B\times B$). 
%For any $1\leq j \leq B$ and any $x\in C$, $\paramstar_{j}(x)$ represents the  $j$-th AP expected signal strength at the position $x$ which is assumed to be independent of $t$. 

The position of the $B$ APs are assumed to be known and denoted by $\{O_{j}\}_{j=1}^{B}$. In order to take into account the perturbations in the signal propagation (due to the fact that radio waves are prone to shadowing, reflections and so on), we propose the following decomposition of $F^{\star}$: for all $x\in\mathcal{C}$ and all $j\in\{1,\cdots,B\}$,
\begin{equation}\label{eqdef:m}
F^{\star}_{j,x} \eqdef \mu^{\star}_{j,x}+\delta^{\star}_{j,x}\eqsp.
\end{equation}
For any $B\times\Rset^{|\mathcal{C}|}$ matrix $A$, we use the shorthand notation $A_j$ for the vector $\{A_{j,x}\}_{x\in\mathcal{C}}$. For any $j\in\{1,\cdots,B\}$, $\mu^{\star}_{j}$ is the average indoor propagation and is such that for all $x\in\mathcal{C}$, $\mu^{\star}_{j,x}$ only depends on the distance between $x$ and $O_{j}$. In the sequel, we use the so-called Friis transmission equation, see \cite{Friis:2006}, given by,
\begin{equation}
\label{eq:def:mu}
\mu^{\star}_{j,x} \eqdef c^{\star}_{1,j} + c^{\star}_{2,j}\log \left\|x-O_{j}\right\|\eqsp,
\end{equation}
where $c^{\star}_{1,j}$ and $c^{\star}_{2,j}$ are parameters depending on the environment and $\log$ is the logarithm to the base e. 

\begin{figure}
\centering
 \subfloat[$\mu_{j}^{\star}$.]{\label{fig:mu_delta_m:mu}\includegraphics[width=0.4\textwidth]{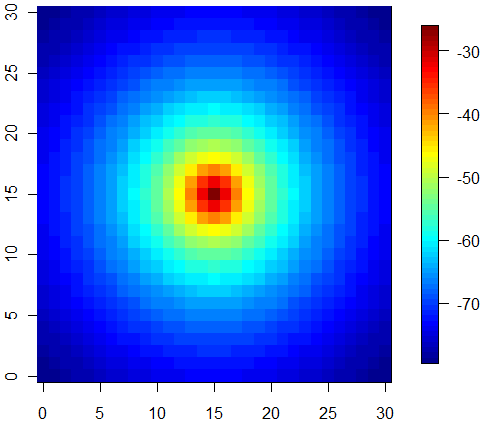}}\\
 \subfloat[$\delta_{j}^{\star}$.]{\label{fig:mu_delta_m:delta}\includegraphics[width=0.4\textwidth]{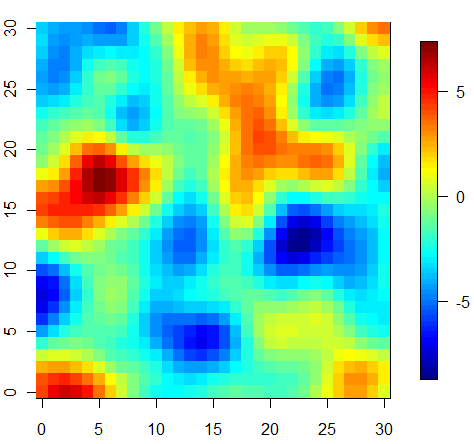}}\\
 \subfloat[$F_{j}^{\star}=\mu_{j}^{\star}+\delta_{j}^{\star}$.]{\label{fig:mu_delta_m:m}\includegraphics[width=0.4\textwidth]{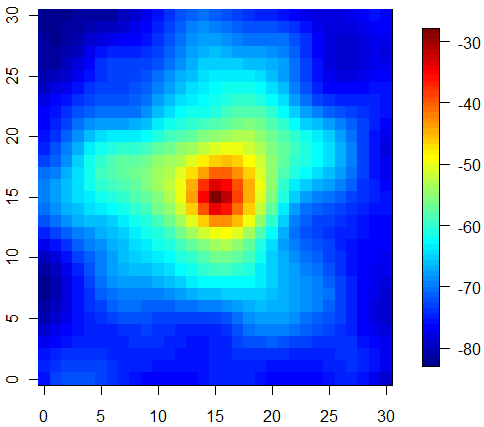}}
\caption{Example of $\delta_{j}^{\star}$ (sampled from $\pi$) and spatial representations of the functions $\mu_{j}^{\star}$ and $F_{j}^{\star}=\mu_{j}^{\star}+\delta_{j}^{\star}$ (in $dBm$)} 
\label{fig:mu_delta_m} 
\end{figure}

$\delta^{\star}_{j}$ is an additive term due to random perturbations such as walls effects (a similar model of WiFi propagation maps using Gaussian processes can be found in \cite{ferris:2006}). It is assumed that the parameters $\{\delta^{\star}_{j}\}_{j=1}^{B}$ are embedded with the prior distribution $\pi$ given, for any $\delta\in\Rset^{B.|\mathcal{C}|}$, by
\[
\pi(\delta)\propto \exp\left\{-\frac{1}{2}\sum_{j=1}^B\delta_j^T\Sigma_j^{-1}\delta_j\right\}\eqsp,
\]
where $\Sigma_{j}$ is assumed to be known and where, for any matrix $A$, $A\trans$ denotes the transpose of A. Figure~\eqref{fig:mu_delta_m} represents $\delta^{\star}_{j}$ (sampled from $\pi$) and the functions  $\mu^{\star}_{j}$ and $F^{\star}_{j}$, defined on the grid $\mathcal{C}=\left\{0,\ldots,30 \right\} \times \left\{0,\ldots,30 \right\}$. The parameters used in this figure are $O_{j} = (15,15)$, and $c^{\star}_{1,j}$, $c^{\star}_{2,j}$ and $\Sigma_{j}$ are given in Section~\ref{sec:simulations}, their values were calibrated after a measurement campaign in an office environment. 

In the sequel, we write $\theta^{\star}\eqdef (c^{\star}_{1}, \ c^{\star}_{2}, \ \delta^{\star}, \ \sigma^{\star,2})$, where $ c^{\star}_{1}\eqdef\{c^{\star}_{1,j}\}_{j=1}^{B}$, $ c^{\star}_{2}\eqdef\{c^{\star}_{2,j}\}_{j=1}^{B}$ and  $ \delta^{\star}\eqdef\{\delta^{\star}_{j}\}_{j=1}^{B}$. For any $x \in \mathcal{C}$, the  distribution of $Y_{t}$ conditionally to $X_{t} =x$ has a density with respect to the Lebesgue measure on $\bigR^{B}$ given, for all $y\eqdef (y_1,\cdots,y_B) \in \bigR^{B}$, by 
\[
g_{\theta^{\star}}(x,y) \eqdef  \dfrac{1}{\sqrt{2\pi\sigma^{\star,2}}^{B}} \prod_{j=1}^{B}\exp\left\{ -\dfrac{1}{2 \sigma^{\star,2}} |y_{j}-F^{\star}_{j,x}|^{2}\right\}\eqsp.
\]

Therefore, $\{X_{t}\}_{t\geq 0}$ is the hidden process of a hidden Markov model observed through the process $\{Y_{t}\}_{t\geq 0}$. The estimation of the mobile device's position $X_{t}$ relies on the knowledge of the map $F^{\star}$. The observations $\{Y_{t}\}_{t \ge 1}$ are used to estimate simultaneously the mobile device's position and $\theta^{\star}$. This simultaneous localization and mapping problem may be seen as an instance of inference in hidden Markov models.
For any positive integer $n$, any observation set $(y_{1},\ldots,y_{n})$, shortly denoted by $y_{1:n}$ and any parameter $\theta = (c_{1},c_{2},\delta,\sigma^{2})$, the likelihood of the observations $L_{\theta}(y_{1:n})$ is given by:
\begin{equation}
L_{\theta}(y_{1:n}) \eqdef \hspace{-.15cm}\sum\limits_{x_{1:n} \in \mathcal{C}^{n}} \nu(x_{1})g_{\theta}(x_{1},y_{1})  \prod\limits_{t=2}^{n} q(x_{t-1},x_{t})g_{\theta}(x_{t},y_{t}) \eqsp.
\end{equation}
Let $n$ be a positive integer and $Y_{1:n}$ be a set of observations, we set as the estimator of $\theta^{\star}$, the maximum a posteriori estimator defined as $\mathrm{argmax}_\theta\, n^{-1} \ell_{\theta} (Y_{1:n})$, where:
\begin{equation}\label{sig:eq:penlikelihood}
\ell_{\theta} (Y_{1:n}) \eqdef \log L_{\theta}(Y_{1:n}) + \log\pi(\delta) \eqsp.
\end{equation}
%where the penalty $ \mbox{pen}(n,\delta) $ on the parameter $\delta$, defined as 
%the logarithm of the Gaussian prior with covariance matrices $\left\{\Sigma_{j} \right\}_{j=1}^{B} $:
%\begin{equation}\label{sig:eq:pen}
 %\mbox{pen}(n,\delta) \eqdef  - \frac{1}{2n} \sum\limits_{j=1}^{B} \delta_{j}^{T} \Sigma_{j}^{-1} \delta_{j}\eqsp, 
%\end{equation}
%comes from the Gaussian prior on the perturbations vectors $\delta_{j}$.
The next section provides a description of the EM algorithm and of online EM algorithms for the computation of maximum likelihood estimators (without the penalty term). In Section~\ref{sec:algo}, we explain how such techniques can be used in our framework.
\section{Online EM}
\label{sec:Online EM}

The EM algorithm is a well-known iterative algorithm to perform maximum likelihood estimation in hidden Markov models~\cite{dempster:laird:rubin:1977}.  Each iteration of this algorithm consists in a E-step where the
expectation of the complete data log-likelihood  (log of the joint distribution of the states and the observations) conditionally to the observations is computed; and a M-step, which updates the parameter estimate. 
%Except for simple models the E-step is intractable and has to be approximated e.g.  by Monte Carlo methods  (see e.g. \cite{cappe:moulines:ryden:2005,fort:moulines:2003}).

Let $Y_{1:n}$ be a fixed set of observations and $\hat{\theta}$ be the current parameter estimate.
\begin{enumerate}[i)]
\item The E-step consists in evaluating the conditional expectation \begin{equation}
\label{eq:Q} Q_{\hat{\theta}}(Y_{1:n};\theta) = \CExpparam{\hat{\theta}}{}{\frac{1}{n}\log p_{\theta}(X_{1:n},Y_{1:n})}{Y_{1:n}}\eqsp,
\end{equation}
where $\log p_{\theta}(X_{1:n},Y_{1:n})$ is the complete data log-likelihood and  $\CExpparam{\hat{\theta}}{}{\cdot}{Y_{1:n}}$ is the conditional expectation given $Y_{1:n}$ when the parameter's value is $\hat{\theta}$. 
\item The M-step updates the current value $\hat{\theta}$ taking the parameter $\theta$ maximizing \eqref{eq:Q}.
\end{enumerate}
The model presented in Section~\ref{sec:framework} belongs to the curved exponential family: there exist functions $S:\Xset^{2}\times\Yset\to\mathcal{S}\subset \mathbb{R}^{d}$, $\phi:\Theta\to \mathbb{R}$ and $\psi:\Theta\to \mathbb{R}^{d}$ such that
\[
 \log q(x,x') + \log g_\theta(x',y) = \phi(\theta) +
  \pscal{S(x,x,',y)}{\psi(\theta)}\eqsp,
\]
where $\pscal{\cdot}{\cdot}$ denotes the scalar product on $\mathbb{R}^{d}$. 
Moreover, there exists a continuous function $\bar \theta:
  \mathcal{S} \to \Theta$ s.t. for any $s \in \mathcal{S}$,
\[
\bar \theta(s) = \mathrm{argmax}_{\theta \in \Theta} \; \left\{ \phi(\theta) +
  \pscal{s}{\psi(\theta)} \right\} \eqsp.
\]

In this case the intermediate quantity defined by \eqref{eq:Q} can be written
\begin{multline}
\label{eq:Qexp} Q_{\hat{\theta}}(Y_{1:n};\theta) =\phi(\theta) \\
+ \pscal{\CExpparam{\hat{\theta}}{}{\frac{1}{n}\sum_{t=1}^{n}S(X_{t-1},X_{t},Y_{t})}{Y_{1:T}}}{\psi(\theta)}\eqsp.
\end{multline}

Therefore, the E-steps amounts to computing only one conditional expectation  $\CExpparam{\hat{\theta}}{}{\frac{1}{n}\sum_{t=1}^{n}S(X_{t-1},X_{t},Y_{t})}{Y_{1:n}}$ when the current parameter's value is $\hat\theta$. The M-step relies simply on the evaluation of $\bar\theta$ at this conditional expectation.
This two steps process is repeated till convergence.
However, when the observations are obtained sequentially or when the E-step relies on a large data set, the EM algorithm might become impractical.  {\em Online} variants of the EM algorithm have been proposed to obtain parameter estimates each time a new observation is available.  In the case of independent and identically distributed (i.i.d.) observations, \cite{cappe:moulines:2009} proposed the first EM based online algorithm. The E-step amounts to computing intermediate quantities known as {\em sufficient statistics} (see below for a explicit definition) and  \cite{cappe:moulines:2009} proposed to replace these computations by a stochastic approximation step. When both the observations and the states take a finite number of values (resp. when the state-space is finite) an online EM-based algorithm was proposed by~\cite{mongillo:deneve:2008} (resp.  by~\cite{cappe:2011}). These algorithms  combine an online approximation of the filtering distributions of the hidden states and a stochastic approximation step to compute an online approximation of the sufficient statistics.  This has been extended to the case of general state-space models with Sequential Monte Carlo algorithms (see~\cite{cappe:2009},~\cite{delmoral:doucet:singh:2010a} and~\cite{lecorff:fort:moulines:2011}). More recently, \cite{lecorff:fort:2011} proposed a block online algorithm in which the parameter estimate is kept fixed on block of observations. The parameter's update then occurs at the end of each block. 

%The convergence of the online EM algorithm for i.i.d.  observations is addressed in~\cite{cappe:moulines:2009}. The limit points are the stationary points of the Kullback-Leibler divergence between the marginal distribution of the
%observation and the model distribution.  \cite{lecorff:fort:2011} gives the asymptotic behavior of its proposed algorithms with rates of convergence. In the case of general state-space models, some insights on the asymptotic behavior are given in~\cite{cappe:2011}.

 In this paper, we use the online variant of the EM algorithm introduced in \cite{lecorff:fort:2011} to perform the parameter estimation and to solve the localization problem presented above. This algorithm, called the {\em Block Online EM} (BOEM)  algorithm relies on the ability to compute sequentially quantities of the form:
\[
\mathsf{S}_{n,\tau}(\theta)\eqdef\CExpparam{\theta}{}{\frac{1}{\tau}\sum_{t=1}^{\tau}S(X_{t+n-1},X_{t+n},Y_{t+n})}{Y_{n:n+\tau}}\eqsp.
\]
Such quantities are called sufficient statistics. The BOEM algorithm uses a sequence of block-sizes $\{\tau_{k}\}_{k\geq 0}$. Define $T_{0}\eqdef 0$ and, for any $k\geq 0$, $T_{k} \eqdef \sum_{i=1}^{k}\tau_{i}$. Let $k$ be a positive integer. Within each block of observations $Y_{T_{k}+1:T_{k+1}}$, the parameter's value $\widehat{\theta}_{k}$ is kept fixed and the sufficient statistic $\mathsf{S}_{T_{k},\tau_{k+1}}(\widehat{\theta}_{k})$ is computed sequentially. The estimate $\widehat{\theta}_{k+1}$ is computed at the end of the block $Y_{T_{k}+1:T_{k+1}}$ through the evaluation of the function $\bar\theta$. 
%The algorithm proposed in \cite{cappe:2011,delmoral:doucet:singh:2010a} updates the parameter each time a new observation is available. It relies on a stochastic approximation step depending on a sequence of block-size $\{\gamma_{n}\}_{n\geq0}$ to update the sufficient statistic. Typically, classic results on stochastic approximation require that this sequence satisfies
%\[
%\sum_{n\geq 0}\gamma_{n} = +\infty\quad\mbox{and}\quad \sum_{n\geq 0}\gamma_{n}^{2} < +\infty\eqsp.
%\]

%
Unlike the traditional use of the EM algorithm where the conditional expectations are computed using forward-backward techniques,  \cite{cappe:2011}, \cite{delmoral:doucet:singh:2010a} and \cite{lecorff:fort:2011} rely on recursive computations of the conditional expectations. Indeed, let $\phi_{n,\theta}$ denotes the filtering distribution of $X_{n}$ given the observations $Y_{1:n}$ when the parameter's value is $\theta$:
\[
\forall x\in\Xset\eqsp,\eqsp\phi_{n,\theta}(x) = P_{\theta}\left(X_{n}=x\vert Y_{1:n} \right)\eqsp.
\]
Following \cite{cappe:2011,delmoral:doucet:singh:2010a}, defining for all $x\in\Xset$ and all $\theta\in\Theta$, 
\[ 
\rho_{n,\theta}(x) =  \CExpparam{\theta}{}{\frac{1}{n}\sum_{t=1}^{n}S(X_{t-1},X_{t},Y_{t})}{Y_{1:n}, \ X_{n}=x}\eqsp, 
\]
we have,
\[
\CExpparam{\theta}{}{\frac{1}{n}\sum_{t=1}^{n}S(X_{t-1},X_{t},Y_{t})}{Y_{1:n}} = \sum_{x \in C}\phi_{n,\theta}(x) \rho_{n,\theta}(x)\eqsp. 
\]
Proposition 1 of \cite{cappe:2011} illustrates the usefulness of this decomposition.
\begin{pr}[of \cite{cappe:2011}] 
$\,$

\noindent\textbf{Initialisation:} 
\\
For all $x \in \mathcal{C}$ and all $\theta \in \Theta$,
\begin{eqnarray*}
\phi_{1,\theta}(x) &=& \frac{\nu(x)g_{\theta}(x,Y_{1})}{\sum\limits_{x' \in \mathcal{C}}\nu(x')g_{\theta}(x',Y_{1})}\eqsp,\\
\rho_{1,\theta}(x) &=&  0\eqsp.
\end{eqnarray*}
\textbf{Recursion:}
\\
For all $t \ge 2$ and all $x \in \mathcal{C}$,
\begin{eqnarray}
\phi_{t,\theta}(x) &=& \frac{ \sum\limits_{x' \in \mathcal{C}}  \phi_{t-1,\theta}(x') q(x'x) g_{\theta}(x,Y_{t}) }{ \sum\limits_{(x',x'') \in \mathcal{C}^{2}}  \phi_{t-1,\theta}(x') q(x'x'') g_{\theta}(x'',Y_{t})}\eqsp,\\
\rho_{t,\theta}(x) &=& \sum\limits_{x' \in \mathcal{C}}  \bigg\{\frac{1}{t} s(x',x,Y_{t}) +  \label{eq:rho}\\ 
& & \left( 1-\frac{1}{t}\right)\rho_{t-1,\theta}(x')\bigg\}.\frac{\phi_{t-1,\theta}(x')q(x',x)}{\sum\limits_{x'' \in \mathcal{C}}\phi_{t-1,\theta}(x'')q(x'',x)}\eqsp.\nonumber
\end{eqnarray}
\end{pr}
Except in simple models (linear Gaussian models and finite state-space HMM), this algorithm requires forward computations which are not available in closed form and which have to be approximated, e.g. using sequential Monte Carlo methods (see \cite{cappe:2009,delmoral:doucet:singh:2010a}). In this case, $\phi_{t,\theta}$ is approximated by weighted samples $\{\widehat{\xi}_{t}^{p},\widehat{\omega}_{t}^{p}\}_{p=1}^{N}$ such that $\widehat{\phi}_{t,\theta}(x) = \sum\limits_{p=1}^{N}\widehat{\omega}_{t}^{p} \delta_{\xi_{t}^{p}}(x)$. In the sequel, $\{\widehat{\xi}_{t}^{p}\}_{p=1}^{N}$ will be referred to as the particle set at time step $t$. Plugging this approximation in (\ref{eq:rho}) yields: 
\begin{equation}
\rho_{t}^{p} = \sum_{\ell=1}^{N}\widehat{\omega}_{t-1}^{\ell}q(\widehat{\xi}_{t-1}^{\ell},\widehat{\xi}_{t}^{p})
\times\frac{\frac{1}{t}s(\widehat{\xi}_{t-1}^{\ell},\widehat{\xi}_{t}^{p},Y_{t})+(1-\frac{1}{t})\rho_{t-1}^{\ell}}{\sum_{\ell=1}^{N}\widehat{\omega}_{t-1}^{\ell}q(\widehat{\xi}_{t-1}^{\ell},\widehat{\xi}_{t}^{p})}\eqsp, \label{eq:rho:p}
\end{equation}
where $\rho_{t}^{p}$ is the approximation of $\rho_{t}$ evaluated at $\widehat{\xi}_{t}^{p}$.
At each time step, the new population of particles is built from the previous population using Algorithm \ref{alg_Bootstrap_filter_recursion} referred to as the {\em bootstrap filter}, see e.g. \cite{cappe:moulines:ryden:2005}. The Bootstrap filter combines sequential importance sampling and sampling importance resampling steps to produce a set of random particles with associated importance weights. Implementations of such procedures are detailed in \cite{cappe:moulines:ryden:2005,cappe:2001,delmoral:2004,doucet:johansen:2009}. 

\begin{algorithm}[h]
\caption{Bootstrap\_filter\_recursion (BFR)}
\label{alg_Bootstrap_filter_recursion}
\begin{algorithmic}[1]
\REQUIRE $\{ \xi_{t-1}^{\ell} , \omega_{t-1}^{\ell} \}_{\ell=1}^{N}$, $Y_{t}$, $\theta$.
\FOR{$p=1$ to $N$} 
\STATE Draw $I$ in ${1,\dots,N}$ with probabilities proportional to $\{\omega_{t-1}^{\ell}\}_{\ell=1}^{N}$.
\STATE Sample $\xi_{t}^{p}\sim  q(\xi_{t-1}^{I},\cdot)$.
\STATE Set $\omega_{t}^{p}\propto g_{\theta}(\xi_{t}^{p},Y_{t})$.  
\ENDFOR
\RETURN $\{ \xi_{t}^{p} , \omega_{t}^{p} \}_{p=1}^{N}$
\end{algorithmic}
\end{algorithm}

This leads to the Algorithm~\ref{alg_SLAM_indoor_BOEM} presented below. Algorithm~\ref{alg_SLAM_indoor_BOEM} is the adaptation of the BOEM to our model. It recursively updates the parameter $\theta$ at the end of each block.
%Algorithm~\ref{alg_SLAM_indoor_BOEM} relies on a sequence of weighted samples $\{(\xi_{t}^{p},\omega_{t}^{p}), p\in\{1,\dots,N\}\}_{t\geq 1}$, henceforth referred to as weighted particles. These particles are used to approximate the recursions not available in closed form. 
The BOEM proposed in \cite{lecorff:fort:2011} also introduced an averaged estimate based on a weighted mean of all the sufficient statistics computed in the past. It is proved in \cite{lecorff:fort:2011} that this averaged estimator has an optimal rate of convergence. Lines~\ref{eq:averagedstat:beg} to~\ref{eq:averagedstat:end} of Algorithm~\ref{alg_SLAM_indoor_BOEM} computes this averaged sufficient statistics and line~\ref{eq:averaged} computes the sequence $\{\widetilde\theta_{k}\}_{k\geq 0}$ of map estimates based on the averaged statistics. The BOEM algorithm is adapted by introducing a second particle system $\{\widetilde{\xi}_{t}^{p},\widetilde{\omega}_{t}^{p} \}_{p=1}^{N}$. This additional particle system is generated using the averaged parameter estimate. As this estimate is supposed to be more accurate than the original estimator  computed on each block, we use the second system of particles to build a better estimator of the device's position. At each time step, we then compute two estimators of the device's position, one for each particle system. Both of them are set as the particle with the greatest importance weight. Line~\ref{eq:suffstat} performs the update of the sufficient statistics and line~\ref{eq:updateparam} the parameter's update at the end of the block.

Finally, in lines~\ref{stabstep:beg} to~\ref{stabstep:end} of Algorithm~\ref{alg_SLAM_indoor_BOEM}, we add a stabilization step (which is not in the original BOEM) which only consists in regularly replace the original map estimate by the averaged one. This step is needed to ensure the convergence as detailed in Section~\ref{sec:simulations}.
\begin{algorithm}[h]
\caption{BOEM\_SLAM\_indoor}
\label{alg_SLAM_indoor_BOEM}
\begin{algorithmic}[1]
\REQUIRE $\theta^{0}$, $\{\tau_{k}\}_{k\geq 1}$, $\{Y_{t}\}_{t\geq 0}$, $N$, $N_{b}\eqsp.$
\STATE Set $\widehat{\theta} =\widetilde{\theta}  = \theta^{0}\eqsp.$
\STATE Sample $\{\widehat{\xi}_{0}^{p}\}_{p = 1}^{N}$ and $\{\widetilde{\xi}_{0}^{p}\}_{p = 1}^{N}$ independently and uniformly in $\mathcal{C}\eqsp.$
\STATE Set $\widehat{\omega}_{0}^{p}= \widetilde{\omega}_{0}^{p}= \frac{1}{N}$  for all $p\in\{1,\dots,N\}\eqsp.$
%\STATE Set $N_{0,2} = N_{2}$
\STATE Set $\rho_{0}^{p} = 0$  for all $p\in\{1,\dots,N\}$, $k = 1$, $T_0=0$, $T_1=\tau_1$.
\medskip
\FOR{all $t\geq 1$} 
\medskip
\STATE {\em Selection and propagation step.}
\STATE Set $\{\widehat{\xi}_{t}^{p},\widehat{\omega}_{t}^{p} \}_{p=1}^{N} = BFR\left(\{\widehat{\xi}_{t-1}^{\ell},\widehat{\omega}_{t-1}^{\ell} \}_{\ell=1}^{N},Y_{t}, \widehat{\theta} \right)$
\STATE Set $\{\widetilde{\xi}_{t}^{p},\widetilde{\omega}_{t}^{p} \}_{p=1}^{N} = BFR\left(\{\widetilde{\xi}_{t-1}^{\ell},\widetilde{\omega}_{t-1}^{\ell} \}_{\ell=1}^{N},Y_{t}, \widetilde{\theta} \right)$
%\STATE Sample $\mathcal{I}$ following a multinomial $\mathcal{M}(N_{2}, \{\tilde{\omega}_{t-1}^{\ell}\}_{\ell=1}^{N_{t-1,2}} )$
%\STATE Set $\tilde{\xi}_{t} = \{x \in C \ \vert \ \exists  \tilde{\xi} \in \mathcal{I} \ s.t. \ \vert x-  \tilde{\xi}\vert \le d_{max}\}$.
%\STATE Set $N_{t,2} = \sharp \tilde{\xi}_{t}$ and $\tilde{\xi}_{t} = \{\tilde{\xi}_{t}^{p}\}_{p=1}^{N_{t,2}}$.
%\FOR{$p=1$ to $N_{t,2}$} 
%\STATE Set $\tilde{\omega}_{t}^{p}\propto g_{\tilde{\theta}}(\tilde{\xi}_{t}^{p},Y_{t})\prod\limits_{\ell \vert \ \vert\tilde{\xi}_{t-1}^{l} -\tilde{\xi}_{t}^{p}\vert \le d_{max}}  q(\tilde{\xi}_{t-1}^{l},\tilde{\xi}_{t}^{p}) \tilde{\omega}_{t-1}^{\ell}$.  
%\ENDFOR
\medskip
\STATE {\em Position estimations.}
\STATE \label{loc:nonaverage}Set $\widehat{p} = \underset{p\in\{1,\dots,N\}}{\argmax}\widehat{\omega}_{t}^{p}$ and  $\widehat{X}_{t} = \widehat{\xi}_{t}^{\widehat{p}}\eqsp.$
\STATE \label{loc:average}Set $\widetilde{p} = \underset{p\in\{1,\dots,N\}}{\argmax}\widetilde{\omega}_{t}^{p}$ and  $\widetilde{X}_{t} = \widetilde{\xi}_{t}^{\widetilde{p}}\eqsp.$
\medskip
\STATE {\em Forward computation of the intermediate quantity.}
\FOR{$p=1$ to $N$} 
\STATE Compute $\{\rho_{t}^{p}\}_{p=1}^{N}$ following (\ref{eq:rho:p}).
\ENDFOR
\medskip
\STATE {\em Map estimation.}
\IF{ $t = T_{k} $}
\STATE Set
\[
\widehat{\mathcal{S}}_{k} =\sum_{p=1}^{N}\widehat{\omega}_{t}^{p}\rho_{t}^{p}\eqsp.
\]\label{eq:suffstat}
\STATE $\widehat{\theta} = \bar{\theta}(\widehat{\mathcal{S}}_{k},\tau_{k})\eqsp.$\label{eq:updateparam}
\STATE Set $\rho_{t}^{p} = 0$  for all $p\in\{1,\dots,N\}\eqsp.$
\IF{$k=1$}\label{eq:averagedstat:beg}
\STATE Set $\widetilde{\mathcal{S}}_{k} = \widehat{\mathcal{S}}_{k}\eqsp.$
\ELSE 
\STATE Set
\[
\widetilde{\mathcal{S}}_{k} = \left(T_{k-1} \widetilde{\mathcal{S}}_{k-1} +\tau_{k} \widehat{\mathcal{S}}_{k} \right)/T_{k} \eqsp.
\]
\ENDIF\label{eq:averagedstat:end}
\STATE $\widetilde{\theta} = \bar{\theta}(\widetilde{\mathcal{S}}_{k},T_{k})\eqsp.$\label{eq:averaged}
\medskip
\STATE {\em Stabilization step.}
\IF{ $ k = 0 \ \text{mod} \ N_{b}$}\label{stabstep:beg}
\STATE Set $\widehat{\theta} = \widetilde{\theta}$
\ENDIF \label{stabstep:end}
\STATE $k=k+1$ and $T_k = T_{k-1}+\tau_k$.
\ENDIF
\medskip
\ENDFOR
\end{algorithmic}
\end{algorithm}

\section{Application of the algorithms to the SLAM in wireless networks}
\label{sec:algo}

In our framework, the objective is the maximisation of the penalized loglikelihood \eqref{sig:eq:penlikelihood}. This task can be performed using a similar technique as the one described in Section~\ref{sec:Online EM} since the additional penalty term only appears in the definition of the function $\bar{\theta}$. Define, for any $(x,y)\in \mathcal{C}\times\bigR^{B} $ and any $j\in \left\{1,\ldots,B\right\}$,
\begin{align*}
s_{1}(x) &\eqdef \{1_{x'}(x)\}_{x'\in \mathcal{C}}\eqsp,\\
s_{2,j}(x,y) &\eqdef \{1_{x'}(x)y_{j}\}_{x'\in \mathcal{C}}\eqsp,\\
s_{3,j}(y) &\eqdef y_{j}^{2}\eqsp.
\end{align*}
The constant $a$ being known, and since our model belongs to the curved exponential family, the penalized intermediate quantity can be written, up to an additive constant, as: 
\begin{multline}
\label{eq:Qexp:pen} Q_{\hat{\theta}}(Y_{1:n};\theta) = -\frac{1}{2n}\sum_{j=1}^{B}\delta_{j}\trans\Sigma_{j}^{-1}\delta_{j}-\frac{B}{2}\log \sigma^{2}
 \\
 -\sum_{j=1}^{B}\frac{\left\{\mathsf{S}_{3,j}-2\pscal{\mathsf{S}_{2,j}}{F_{j}}+\pscal{\mathsf{S}_{1}}{F_{j}^{2}}\right\}}{2\sigma^{2}}\eqsp,
\end{multline}

%Where $\frac{1}{T}\sum_{j=1}^{B}\delta_{j}\trans\Sigma_{j}^{-1}\delta_{j}$ is the Gaussian prior on $\{\delta_{j}\}_{j=1}^{B}$ presented in the preceding section, 
where, $F_j^2\eqdef \{F_{j,x}^2\}_{x\in\mathcal{C}}$ and 
\begin{equation*}
\mathsf{S}_{1} \eqdef \frac{1}{n}\CE[Y_{1:n}]{}{\widehat{\theta}}{\sum_{t=1}^{n}s_{1}(X_{t})}\eqsp,
\end{equation*}
and, for all $j\in\{1,\dots,B\}$,
\begin{align*}
\mathsf{S}_{2,j} &\eqdef \CE[Y_{1:n}]{}{\widehat{\theta}}{\frac{1}{n}\sum_{t=1}^{n}s_{2,j}(X_{t},Y_{t})}\eqsp,\\
\mathsf{S}_{3,j} &\eqdef \frac{1}{n}\sum_{t=1}^{n}s_{3,j}(Y_{t})\eqsp.
\end{align*}

For any $\mathsf{S}=\left(\mathsf{S}_{1},\ \{\mathsf{S}_{2,j}\}_{j=1}^{B},\ \{\mathsf{S}_{3,j}\}_{j=1}^{B}\right)\in \left[0,1 \right]^{|\mathcal{C}|} \times \bigR^{\vert \mathcal{C}\vert\times B } \times \left(\bigR^{\star}_{+}\right)^{\vert \mathcal{C}\vert\times B }$, we denote by $\bar{\theta}(\mathsf{S},n)$ one of the parameter $\theta=( c_{1}, c_{2},\delta, \sigma^{2})$ maximizing the expression:
\begin{multline*}
-\sum_{j=1}^{B}\frac{\left\{\mathsf{S}_{3,j}-2\pscal{\mathsf{S}_{2,j}}{F_{j}}+\pscal{\mathsf{S}_{1}}{F_{j}^{2}}\right\}}{2\sigma^{2}}\\
+\frac{1}{n}\log\pi(\delta) -\frac{B}{2}\log \sigma^{2}\eqsp.
\end{multline*}
 
As mentioned in Section~\ref{sec:framework}, for any $j\in\{1,\dots,B\}$, $ F_{j}$ is written as $ F_{j} \eqdef \mu_{j}+\delta_{j}$. In these experiments, the whole set of parameters ($c_{1}^{\star}$, $c_{2}^{\star}$ and $\sigma^{\star,2}$) and the unknown perturbation Gaussian fields $\{\delta_{j}^{\star}\}_{j=1}^{B}$ are estimated using Algorithm~\ref{alg_SLAM_indoor_BOEM}. For all $j \in \{1,\dots, B\}$, we write $D_j\eqdef\{\log\|x-O_j\|\}_{x\in\mathcal{C}}$ and
\begin{eqnarray*}
M_{0,j} &\eqdef& \left[ \text{diag}(\mathsf{S}_{1}) + \frac{\sigma^{2}}{n+1} \Sigma_{j}^{-1}\right]\eqsp,\label{etape1}\\
M_{1,j} &\eqdef& \text{diag}(\mathsf{S}_{1}) \left[I-M_{0,j}^{-1}\text{diag}(\mathsf{S}_{1}) \right]\eqsp,\label{etape2}\\
M_{2,j} &\eqdef& I - \text{diag}(\mathsf{S}_{1}) M_{0,j}^{-1}\eqsp,\label{etape3}\\
W_{1,j} &\eqdef& \mathds{1}\trans M_{1,j}\mathds{1}\eqsp,\label{etape4} \\
W_{2,j} &\eqdef& \mathds{1}\trans M_{1,j} D_{j}\eqsp,\label{etape5}\\
%W_{3,j} &\eqdef& W_{2,j}\label{etape6}\eqsp,\\
W_{3,j} &\eqdef& D_{j}\trans M_{1,j} D_{j}\eqsp,\label{etape7}\\
d_{j}&\eqdef& W_{1,j}W_{4,j} - W_{2,j}^2\eqsp.
\end{eqnarray*}

Thus, $\theta =  \bar{\theta}(\mathsf{S},n)$ is given, by $\theta=( c_{1},c_{2},\delta, \sigma^{2})$ where, for all $j \in \{1,\dots, B\}$,
\begin{eqnarray*}
 c_{1,j} &=& d_{j}^{-1} \left[W_{3,j} \mathds{1}\trans - W_{2,j} D_{j}\trans \right] M_{2,j} \mathsf{S}_{2,j}\label{etape8}\eqsp,\\
 c_{2,j} &=& d_{j}^{-1}\left[-W_{2,j} \mathds{1}\trans + W_{1,j} D_{j}\trans \right] M_{2} \mathsf{S}_{2,j}\label{etape9}\eqsp,\\
\delta_{j} &=& M_{0,j} \left[\mathsf{S}_{2,j} - \text{diag}(\mathsf{S}_{1}) (c_{1,j} \mathds{1} +  c_{2,j}D_{j}) \right]\label{etape10}\eqsp,\\
 F_{j} &=&   c_{1,j}\mathds{1} +  c_{2,j}D_{j} +\delta_{j}\label{etape11}
\end{eqnarray*}
and 
\begin{equation*}
\sigma^{2}= \frac{1}{B} \sum\limits_{j=1}^{B} \left\{F_{j}^{\trans} \text{diag}(\mathsf{S}_{1})  F_{j} - 2\mathsf{S}_{2,j}\trans   F_{j} + \mathsf{S}_{3,j} \right\}\eqsp.\label{etape13}
\end{equation*}

\section{Experiments}\label{sec:simulations}
\subsection{Simulated data}
In this section, the performance of the proposed BOEM algorithm is illustrated with simulated data.
All experiments are performed on the grid $\mathcal{C} = \{0,\dots,30\}\times \{0,\dots,30\}$. We use $B = 17$ APs, each AP being modelled by the same coefficients $c_{1}^{\star}$ and $c_{2}^{\star}$, see \eqref{eq:def:mu},
\[
\forall j\in \{1,\dots,B\},\eqsp c_{1,j}^{\star} = -26\quad\mbox{and}\quad c_{2,j}^{\star} = -17.5\eqsp.
\] 
For all $j\in \{1,\dots,B\}$, $\Sigma_{j}$ is a Gaussian covariance function defined by $\Sigma_{j}(x,x')\eqdef v_{1}*\mbox{exp}(-|x-x'|^2/v_{2})$ with $v_{1} =  10$ and $ v_{2} =  18$. The variance of the observation noise is $\sigma^{\star,2} = 25$. The variance of the transition kernel defined in \eqref{eq:def:q} is chosen such that $a=6$.

\medskip

All runs are started with the same initial estimates $\theta^{0} =\left(c^{0}_{1},c^{0}_{2},\delta^{0},\sigma^{0,2}\right) $ where, $\delta^{0}=0$
\[
\forall j\in\{1,\dots,B\},\eqsp  c^{0}_{1,j} =-10\eqsp,\eqsp  c^{0}_{2,j} =-30\quad\mbox{and}\quad  \sigma^{0,2}=30\eqsp.
\]

\begin{figure}
\centering
 \subfloat[$0.8$-quantile of the distance between the true localization and the estimated position. The localization error is given with the nonaveraged estimate (dotted line), the averaged estimate (dashed line) and the optimal estimate (bold line).]{\label{fig:map:loc:withoutinj:loc}\includegraphics[width=0.5\textwidth]{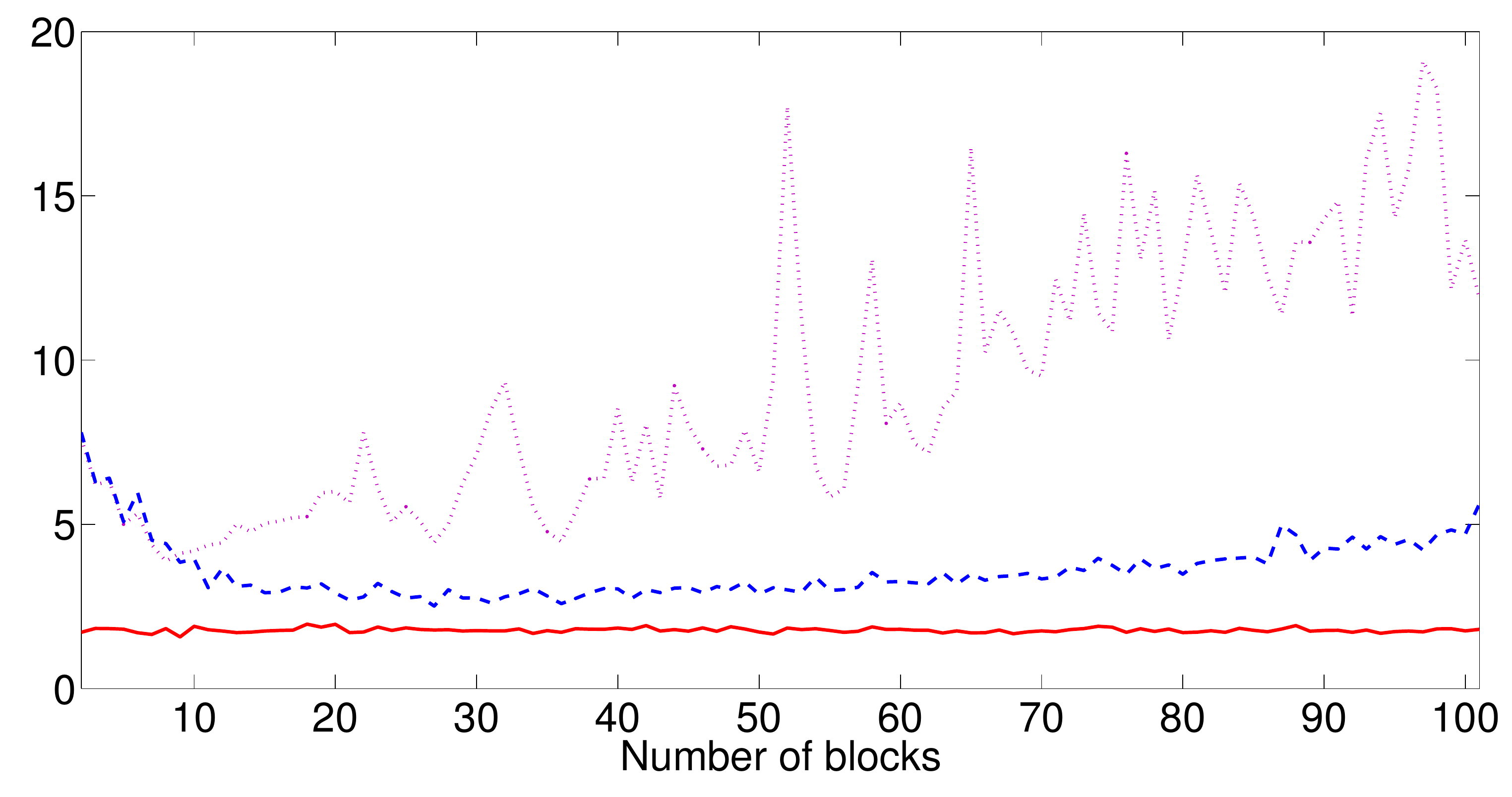}}\\
 \subfloat[Mean $\mathrm{L}_{1}$ error on the map estimate with the initial estimate (dotted line) and the averaged estimate (dashed line).]{\label{fig:map:loc:withoutinj:map}\includegraphics[width=0.5\textwidth]{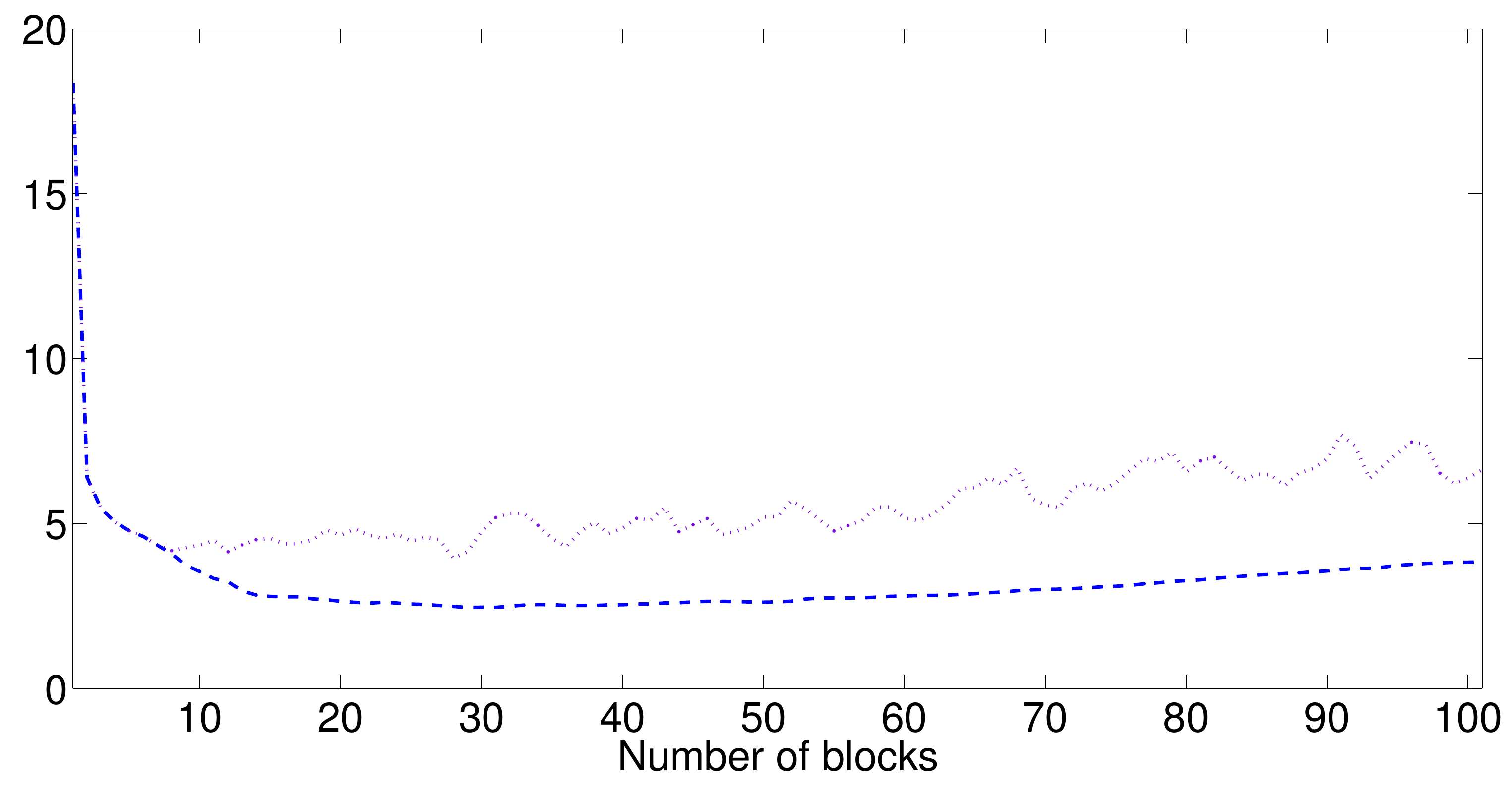}}
\caption{Errors on the map estimation and localization processes with the original algorithm.} 
\label{fig:map:loc:withoutinj} 
\end{figure}

The number of particles $N=25$ is kept fixed and the initial position of each particle is chosen randomly and uniformly in $\mathcal{C}$. For each map $F_{j}^{\star}$, the estimation error is set as the normalized $\mathrm{L}_{1}$ error, such that the distance of a given map $F_{j}$ from the true map $\paramstar_{j}$ is
\[
\epsilon_{j} \eqdef \frac{1}{|\mathcal{C}|}\sum_{x\in \mathcal{C}}\left|F_{j,x}-F^{\star}_{j,x}\right|\eqsp,
\]
and the error displayed is the mean over all maps:
\[
\bar\epsilon \eqdef \frac{1}{B}\sum_{j=1}^{B}\epsilon_{j}\eqsp,
\]
The block sizes are given by
\[
\forall k\in\Nset,\eqsp \tau_{k} = 10k+500\eqsp.
\]
On each block, the localization error is set as the $0.8$-quantile of the distance between the true localization and the estimated position. 
Figure~\ref{fig:map:loc:withoutinj} displays the error on the estimation of the maps and on the localization when the stabilization step in Algorithm \ref{alg_SLAM_indoor_BOEM} is omitted (lines 29 to 31). This case corresponds to the BOEM algorithm. The localization part is dealt with using two different  procedures.
\begin{itemize}
\item {\em Nonaveraged estimate}: the estimate is given with the original particle system (see line~\ref{loc:nonaverage} of Algorithm~\ref{alg_SLAM_indoor_BOEM}).
\item {\em Averaged estimate}: the estimate is given with a second particle system run with the average estimation of the map (see line~\ref{loc:average} of Algorithm~\ref{alg_SLAM_indoor_BOEM}).
\end{itemize}
In order to give fair results, the optimal estimate is shown, i.e. the estimated position given with a particle system run with the true maps $F_{j}^{\star}$, $j\in\{1,\dots,B\}$. 

\begin{figure}
\centering
 \subfloat[$0.8$-quantile of the distance between the true localization and the estimated position with the stabilization process. The localization error is given with the nonaveraged estimate (dotted line), the averaged estimate (dashed line) and the optimal estimate (bold line).]{\label{fig:map:loc:withinj:loc}\includegraphics[width=0.5\textwidth]{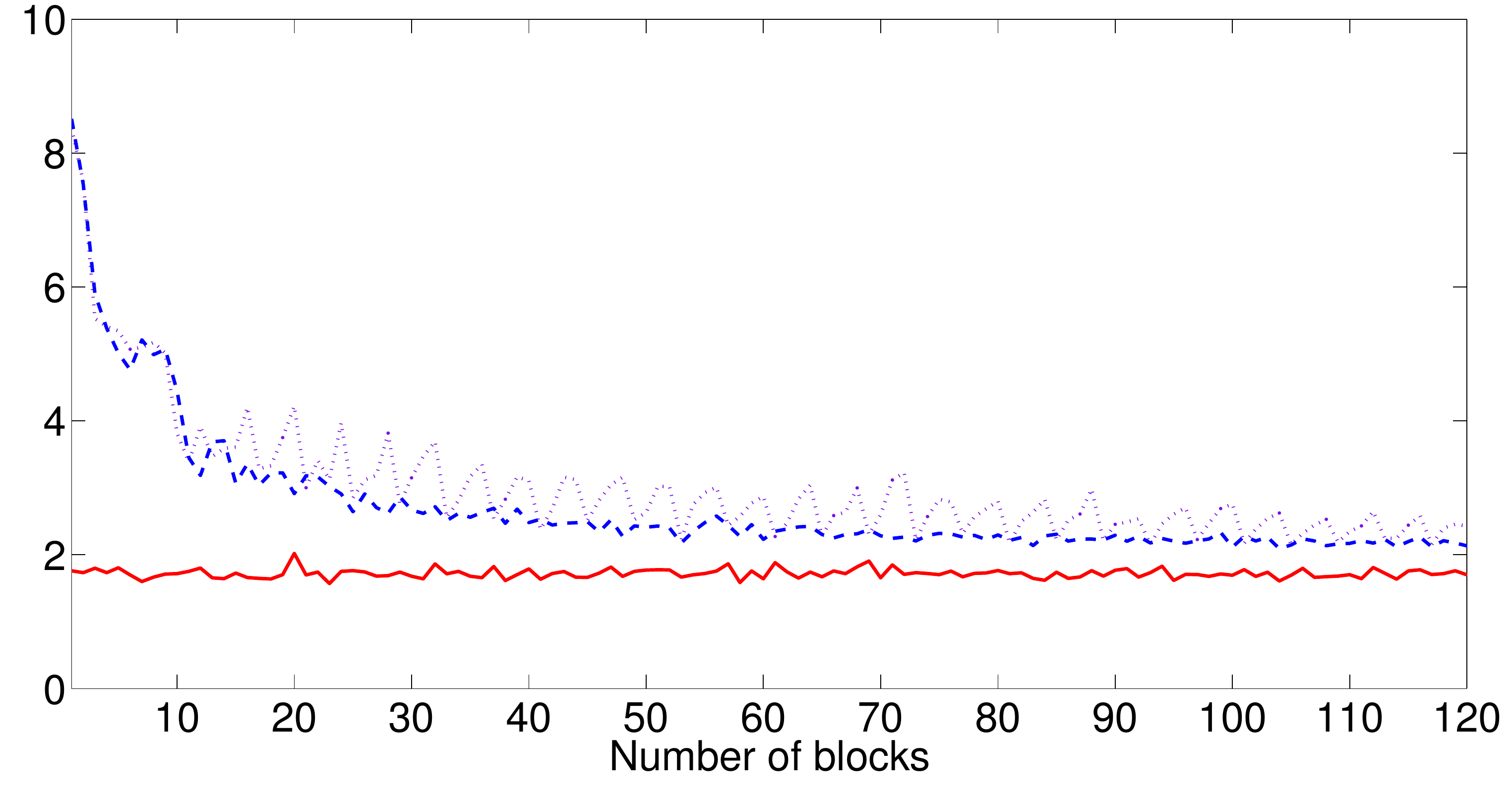}}\\
 \subfloat[Mean $\mathrm{L}_{1}$ error on the map estimate with the initial estimate (dotted line) and the averaged estimate (dashed line) with the stabilization process.]{\label{fig:map:loc:withinj:map}\includegraphics[width=0.5\textwidth]{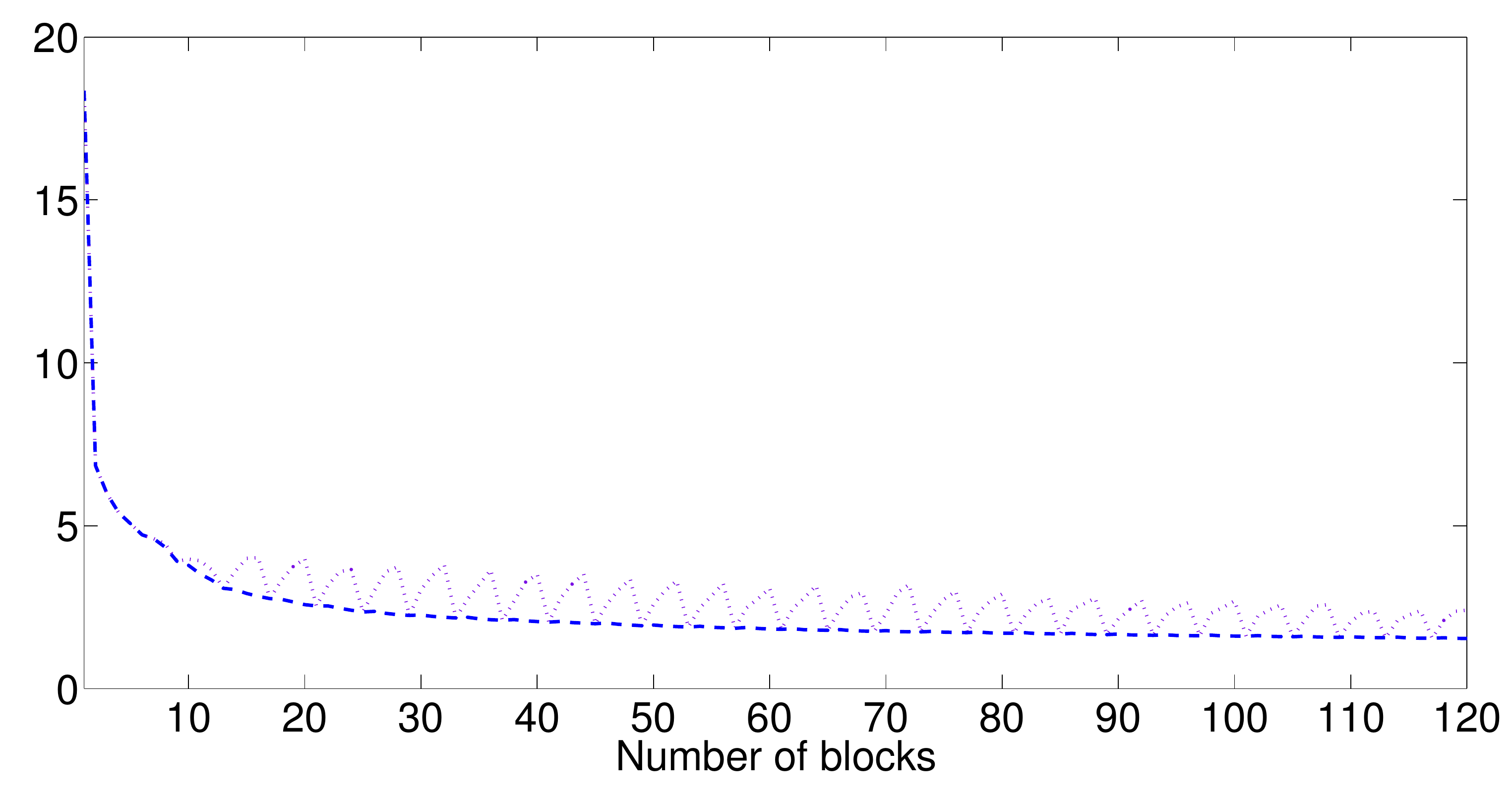}}
\caption{Errors on the map estimation and localization procedure with the stabilized algorithm.} 
\label{fig:map:loc:withinj} 
\end{figure}

As shown in Figure~\ref{fig:map:loc:withoutinj:loc} the estimated position does not converge as the number of blocks (i.e. as the number of estimations) increases. After $50$ blocks (about $40000$ observations) the position, which is badly estimated, does not provide good map estimates which increases the error on the averaged map estimate. Figure~\ref{fig:map:loc:withoutinj:map} displays the error on the map estimate. It is clear that both the estimate and its averaged version do not converge. This convergence problem of the BOEM algorithm can be due to the curse of dimensionality that can occur when the number of parameters to estimate is high. Moreover, the higher the parameter space dimension is, the more likely EM based algorithms are prone to converge towards local minima (see \cite{cappe:moulines:ryden:2005}). To overcome this difficulty, we propose to use the good behaviour of the averaged map estimate during the first $50$ blocks. The map estimate is regularly replaced by its averaged version, see lines~\ref{stabstep:beg} to~\ref{stabstep:end} of Algorithm~\ref{alg_SLAM_indoor_BOEM}. This will prevent the map estimate from diverging and thus, this will reduce the error on the estimated position. In Figure~\ref{fig:map:loc:withinj}, this stabilization process is performed each time $N_{b}=5$ blocks have been used. As shown by Figure~\ref{fig:map:loc:withinj:loc} and Figure~\ref{fig:map:loc:withinj:map}, this greatly improves the performance of the estimation of both the maps and the localization. Hence, the proposed algorithm is based on this stabilization procedure and uses the averaged position estimate to perform the localization part.

\begin{figure}
\centering
\includegraphics[width=0.5\textwidth]{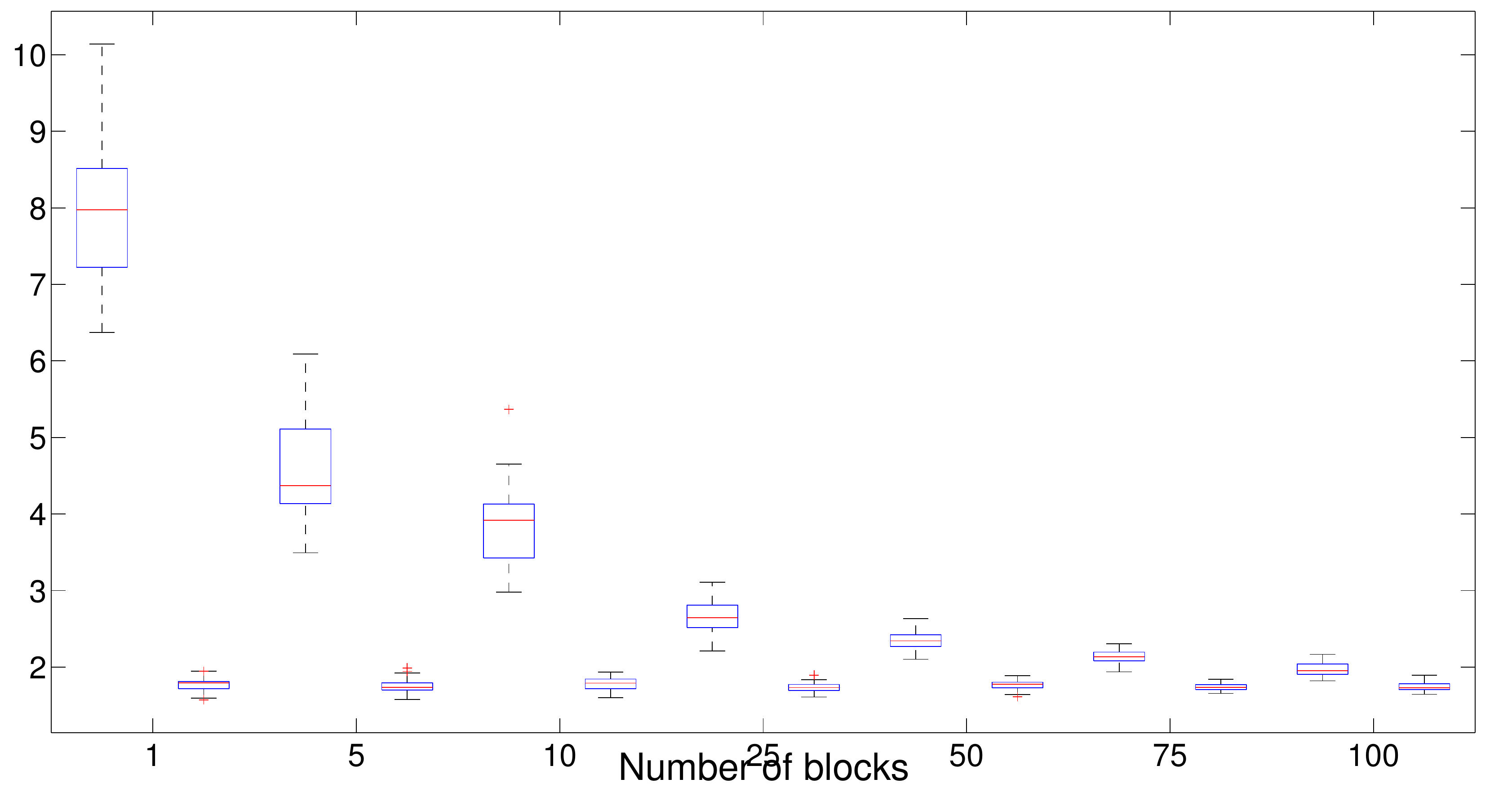}
\caption{Boxplots of the localization error given by the stabilized algorithm with the averaged estimate (left) and the optimal estimate (right) as a function of the number of blocks.}
\label{fig:loc_stab}
\end{figure}

\begin{figure}
\centering
\includegraphics[width=0.5\textwidth]{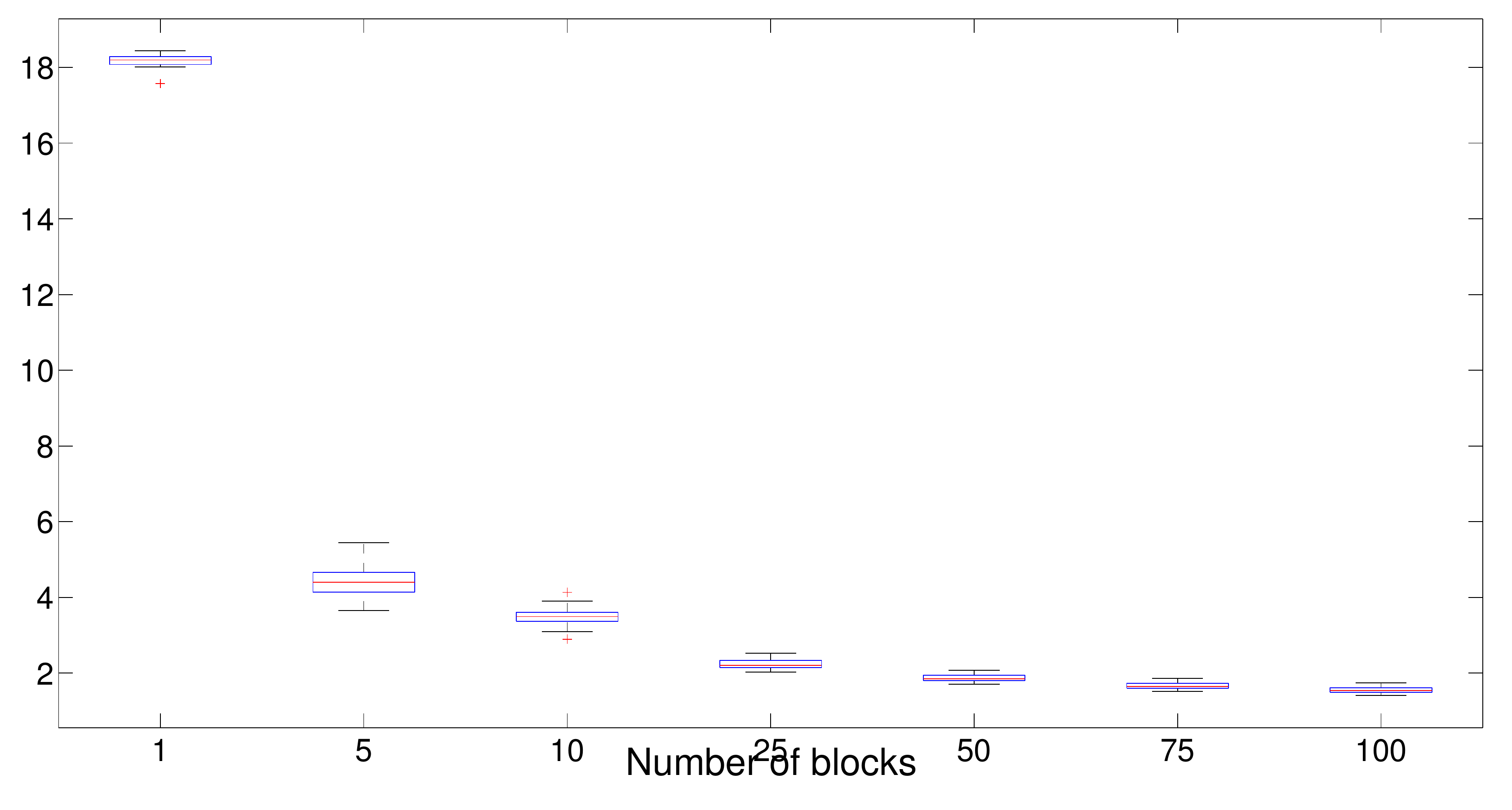}
\caption{Boxplots of the mean $\mathrm{L}_{1}$ error on the map estimate with the stabilized algorithm and the averaged estimate as a function of the number of blocks.}
\label{fig:map_stab}
\end{figure}

Figures~\ref{fig:loc_stab} and \ref{fig:map_stab} illustrate the performance of the algorithm  for the localization and for the estimation of the maps over $50$ independent Monte Carlo runs. In Figure~\ref{fig:loc_stab}, the optimal localization error (i.e. when the maps are known) is also displayed. The convergence of the localization error to the optimal error is almost reached after $100$ blocks (about $100000$ observations). Similarly, the error for the estimation of the maps given by the averaged algorithm goes on decreasing after $100$ blocks (the decrease is slower after $75$ blocks).  

\subsection{True data}\label{sig:sec:truedata}
\begin{figure}
\centering
\includegraphics[width=0.5\textwidth]{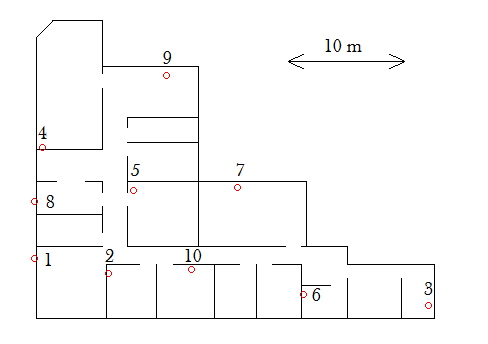}
\caption{Map of the indoor environment used for the test with the position the access points (red circles) and their associated identification numbers.}
\label{fig:map:idservices}
\end{figure}

\begin{figure}
\centering
\includegraphics[width=0.5\textwidth]{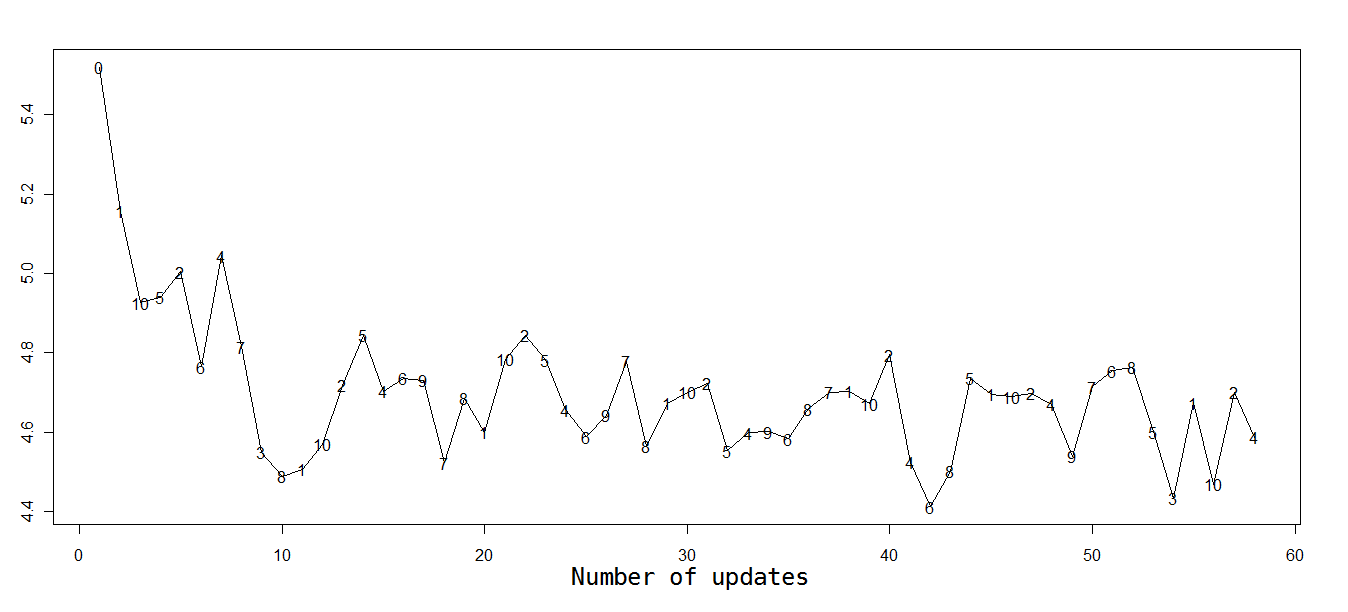}
\caption{0.8-quantile of the distance (in meter) between the true localization and the averaged estimate obtained with the stabilized algorithm. The localization error is computed on the test data sample each time one of the estimated map is updated.}
\label{fig:map:expresults}
\end{figure}

In this section, the behaviour of our SLAM algorithm is illustrated in a real situation. 10 access points are set up in an office environment (Figure~ \ref{fig:map:idservices} represents a map of this environment as well as the position of the access points). The map is discretized using a grid $\mathcal{C} \subset [0,30]\times [0,30]$. 
The variance $\sigma^{\star,2}$ is assumed to be known and its value ($\sigma^{\star,2} = 25dBm^{2}$) is calibrated using a measurement campaign at a fixed position.
Around $T=20000$ measures of the RSSI have been made on the map using a WiFi device. Algorithm \ref{alg_SLAM_indoor_BOEM} produces position estimates 
%$\{\hat{X}_{t}\}_{t \ge 0}$ and $\{\hat{\tilde{X}}_{t}\}_{t \ge 0}$
but we do not have a direct access to the real position and thus cannot observe the localization error. To overcome this difficulty, a test data sample is built  by producing measures along 12 paths in the environment such that, for each measure, the associated position is registered. The test data sample is made of $T_{\text{test}}=1100$ measures: $\{X_{t}^{\text{test}},Y_{t}^{\text{test}}\}_{t=1}^{T_{\text{test}}}$.
The test data sample is used to compare the localization accuracy provided by different values of the parameter $F$. Note a major difference between the model given in Section~\ref{sec:framework} and the real data situation. For any measure $Y$ sent by the device, only several APs are represented in $Y$. Therefore,  the maps $\widetilde{F}_{j}$, $j\in\{1,\ldots,B\}$ are not estimated simultaneously as, for any time step $t$, two APs might appear a different number of times in $Y_{1:t}$. We thus slightly modify Algorithm~\ref{alg_SLAM_indoor_BOEM} by introducing specific blocks and measure counters relatively to each AP. Each time the value of $\widetilde{F}_{j}$ for any $j \in \{1,\ldots,B\}$ is updated using a block of the $20000$ measures, we submit the new estimator $\widetilde{F}$ to the test data sample:  Algorithm~\ref{alg_SLAM_indoor_BOEM} is run on the test data sample $\{Y_{t}^{\text{test}}\}_{t=1}^{T_{\text{test}}}$. Only the averaged particle system is computed and no parameter update is performed. We can then compute the localization error relatively to the test data sample as the 0.8-quantile of the error between $\{X_{t}^{\text{test}}\}_{t=1}^{T_{\text{test}}}$ and the averaged position estimate. Figure \ref{fig:map:expresults} displays the results of this experiment by representing the 0.8-quantile of the error as a function of the number of updates. The numbers on the graph in Figure \ref{fig:map:expresults} indicate which AP were updated for each update. The initial map estimates are given, for any $j\in\{1,\ldots,10\}$ by $  c^{0}_{1,j} = -26$ and $ c^{0}_{2,j} = -17.5$ and $\delta_{0} = 0$.
Despite the relatively small test sample size, Figure \ref{fig:map:expresults} shows that the localization error seems to adopt the same behaviour as the localization error for the simulated data. The parameter $F_{j}$ was updated a maximum of $7$ times (for AP $j=10$ for instance) and a minimum of $2$ times (for AP $j=3$).
\begin{figure}
\centering
\includegraphics[width=0.5\textwidth]{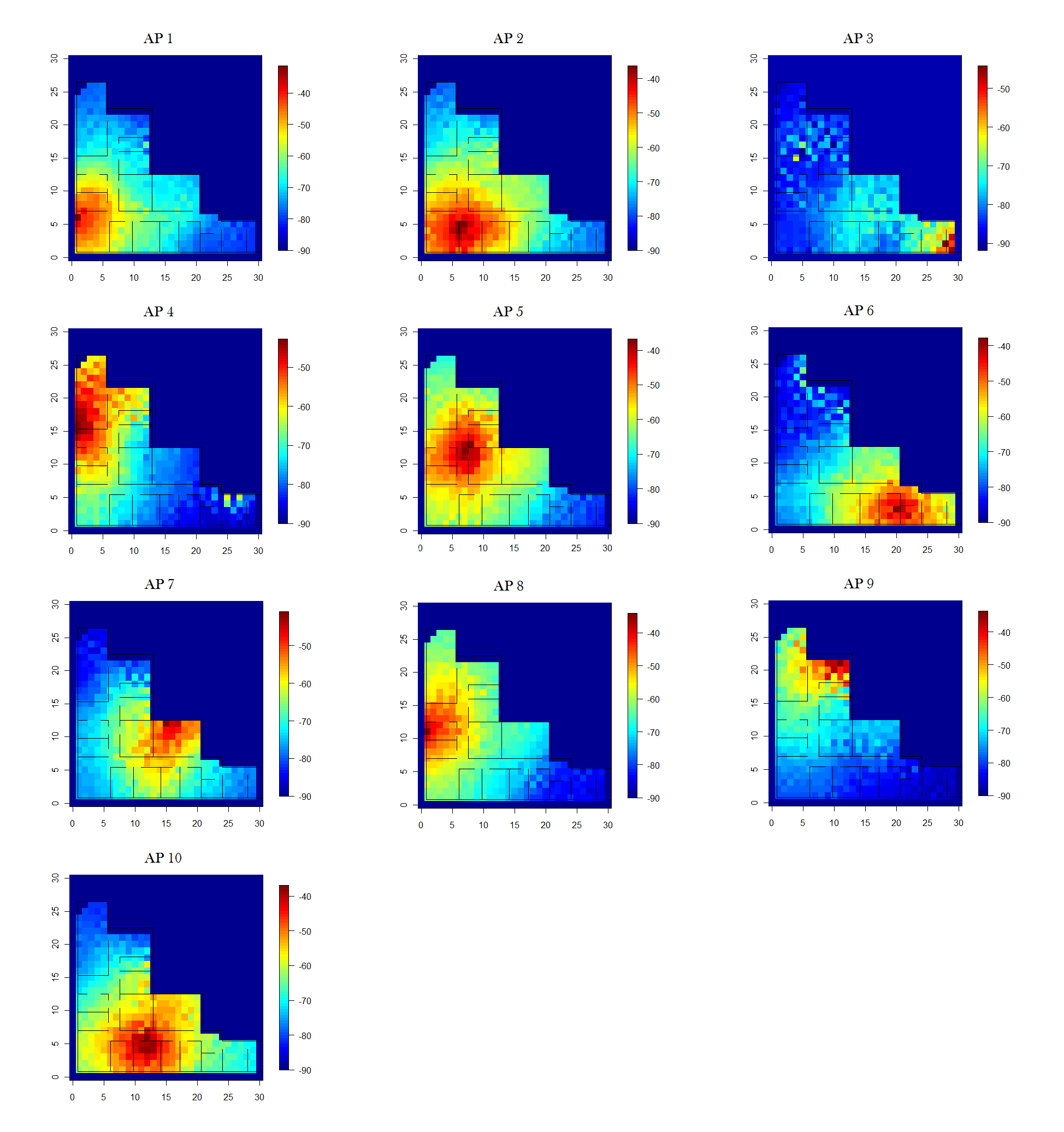}
\caption{Graphical representation of the final propagation maps estimations: $\{\widetilde{F}_{j}\}_{j=1}^{10}$ (in $dBm$). }
\label{fig:map:mesti}
\end{figure}
Figure \ref{fig:map:mesti} represents the final estimate of the propagation maps $\widehat F_{j}$, $j\in\{1,\ldots,10\}$. For some map estimates (see for instance the access points $1$, $4$ and $7$), the signal strength drops when passing walls while the walls are responsible for a part of the indoor waves propagation disturbances.

\section{Conclusion}
In this paper we propose a stabilized version of the BOEM algorithm to estimate the signal propagation maps needed in any WiFi based localization system. The main difference with the existing solutions is that these propagation maps are estimated using the data sent by the mobile device originally used for localization purposes. On the contrary, the existing WiFi based localization systems establish these propagation maps either in a deterministic way or by running a previous hand made survey. In case of environmental modifications, the propagation maps are thereby changed. Our technique can easily be adapted to these changes by regularly reinitializing the sufficient statistics while hand made survey based systems can not take into account these modifications without renewing the survey. However, further tests are needed to evaluate the accuracy provided by our method and to compare it with other methods. Many elements should be analyzed such as the number and the position of the access points, the size of the environment or the materials constituting the obstacles in the environment. 

\bibliographystyle{IEEEtran}
\bibliography{refsSLAMindoor}

\end{document}